\documentclass[11pt]{article}
\usepackage{latexsym}
\usepackage[all]{xy}
\usepackage{amsmath}
\usepackage{amssymb}
\usepackage{amsfonts}

\newcommand  {\dst}   {\mathbf{dSt}_{/k}}

\newtheorem{thm}{Theorem}[section]
\newtheorem{prop}[thm]{Proposition}

\newtheorem{df}[thm]{Definition}

\newtheorem{conj}[thm]{Conjecture}

\begin{document}

\title{\textbf{A note on Chern character, loop spaces and derived algebraic geometry}} 
\bigskip
\bigskip

\author{\bigskip\\
Bertrand To\"en \\
\small{Institut de Math\'ematiques de Toulouse}\\
\small{UMR CNRS 5219} \\
\small{Universit\'{e} Paul Sabatier, Bat. 1R2 Toulouse Cedex 9}\\
\small{France}\\
\bigskip \and
\bigskip \\
Gabriele Vezzosi \\
\small{Dipartimento di Matematica Applicata}\\
\small{Universit\`a di Firenze}\\
\small{Italy}\\
\bigskip}

\date{March 2008}

\maketitle

\begin{abstract}In this note we present a work in progress whose main purpose is to 
establish a categorified version of sheaf theory. We present
a notion of \emph{derived categorical sheaves}, which is a categorified version of 
the notion of complexes of sheaves of $\mathcal{O}$-modules on schemes, as
well as its quasi-coherent and perfect versions. We also explain how ideas from 
derived algebraic geometry and higher category theory can be used in order
to construct a Chern character for these categorical sheaves, which is
a categorified version of the Chern character for perfect complexes with values
in cyclic homology. Our construction uses
in an essential way the \emph{derived loop space} of a scheme $X$, which is a \emph{derived scheme} whose
theory of functions is closely related to cyclic homology of $X$. This work can be seen as an attempt
to define algebraic analogs of elliptic objects and characteristic classes for them. The present text is an overview
of a work in progress and details will appear elsewhere.
\end{abstract}

\tableofcontents

\section{Motivations and objectives}

The purpose of this short note is to present a construction of a Chern character
defined for certain sheaves of categories rather than sheaves of modules (e.g. vector 
bundles or coherent sheaves). This is part of a more ambitious project to 
develop a general theory of categorical sheaves, in the context of algebraic
geometry but also in topology, which is supposed to be a categorification 
of the theory of sheaf of modules. 
Our original motivations for starting such a project
come from elliptic cohomology, which we now explain briefly. \\

\begin{center} \textbf{From elliptic cohomology to categorical sheaves} \end{center}

To any (complex oriented) generalized cohomology theory $E_{*}$ (defined on 
topological spaces) is associated an integer called its
chromatic level, which by definition is the height of the
corresponding formal group. The typical generalized cohomology theory 
of height zero is singular cohomology and is represented by the
Eleinberg-MacLane spectrum $H\mathbb{Z}$. The typical generalized cohomology theory
of height 1 is complex K-theory which is represented by the spectrum $BU\times \mathbb{Z}$. 
A typical cohomology theory of height 2 is
represented by an elliptic spectrum and is called elliptic cohomology. These elliptic cohomologies can be combined altogether 
into a spectrum $\mathsf{tmf}$ of \emph{topological modular forms} (we recommend the excellent survey
\cite{lu} on the subject).
The cohomology theories $H\mathbb{Z}$ and $BU\times \mathbb{Z}$ are rather well understood, in the sense
that for a finite CW complex $X$ is possible to describe the groups $[X,H\mathbb{Z}]$ and
$[X,BU\times \mathbb{Z}]$ easily in terms of the topology of $X$. Indeed, 
$[X,H\mathbb{Z}]\simeq H^{0}(X,\mathbb{Z})$ is the group of continuous functions
$X \longrightarrow \mathbb{Z}$. In the same way, 
$[X,BU\times \mathbb{Z}]=K_{0}^{top}(X)$ is the Grothendieck group of complex
vector bundles on $X$. As far as we know, it is an open question to describe
the group $[X,\mathsf{tmf}]=Ell_{0}(X)$, or the groups $[X,E]$ for some elliptic spectrum $E$,  in similar terms, e.g. as the Grothendieck group of some
kind of geometric objects over $X$ (for some recent works in this direction, see \cite{bdr} and \cite{st-te}). 

It has been observed by several authors that the chromatic level of the cohomology theories
$H\mathbb{Z}$ and $BU\times \mathbb{Z}$ coincide with a certain \emph{categorical level}. More precisely, 
$[X,H\mathbb{Z}]$ is the set of continuous functions $X \longrightarrow \mathbb{Z}$. In this description
$\mathbb{Z}$ is a discrete topological space, or equivalently a set, or equivalently 
a \emph{$0$-category}. In the same way, classes in $[X,BU\times \mathbb{Z}]$ can be represented by 
finite dimensional complex vector bundles on $X$. A finite dimensional complex vector bundle on $X$ is a 
continuous family of finite dimensional complex vector spaces, or equivalently a \emph{continuous map} 
$X \longrightarrow \underline{Vect}$, 
where $\underline{Vect}$ is the \emph{1-category} of finite dimensional complex vector spaces. Such an 
interpretation of vector bundles can be made rigorous if $\underline{Vect}$ is considered
as a topological stack. It is natural to expect that $[X,\mathsf{tmf}]$ is related in one way or another to 
$2$-categories, and that classes in $[X,\mathsf{tmf}]$ should be represented by certain
continuous applications $X \longrightarrow 2-\underline{Vect}$, where 
now $2-\underline{Vect}$ is a $2$-category (or rather a topological $2$-stack). The notation
$2-\underline{Vect}$ suggests here that $2-\underline{Vect}$ is a categorification 
of $\underline{Vect}$, which is itself a categorification of $\mathbb{Z}$ (or rather of
$\mathbb{C}$). If we follow this idea further the typical generalized cohomology theory
$E$ of chromatic level $n$ should itself be related to $n$-categories in the sense that 
classes in $[X,E]$ should be represented by continuous maps $X \longrightarrow n-\underline{Vect}$, 
where $n-\underline{Vect}$ is now a certain topological $n$-stack, which is supposed to be
an $n$-categorification of the $(n-1)$-stack $(n-1)-\underline{Vect}$.

This purely formal observation relating the chromatic level to a certain categorical level
is in fact supported by at least two recent results. On the one hand, J. Rognes 
stated the so-called \emph{red shift conjecture}, which from an intuitive point of view stipulates that 
if a commutative ring spectrum $E$ is of chromatic level $n$ then its $K$-theory spectrum
$K(E)$ is of chromatic level $n+1$ (see \cite{ro}). Some explicit computations
of $K(BU\times \mathbb{Z})$ proves a major case of this conjecture for $n=1$ (see \cite{bdr}). Moreover, 
$K(BU\times \mathbb{Z})$ can be seen to be the $K$-theory spectrum of the $2$-category
of complex $2$-vector spaces (in the sense of Kapranov-Voevodsky). This clearly shows
the existence of an interesting relation between elliptic cohomology and 
the notion of \emph{$2$-vector bundles} (parametrized version of the
notion $2$-vector spaces), even though the precise relation remains unclear at the moment.
On the other hand, the fact that topological $K$-theory is obtained as the Grothendieck group
of vector bundles implies the existence of equivariant $K$-theory by using 
equivariant vector bundles. It is important to notice here that the spectrum
$BU\times \mathbb{Z}$ alone is not enough to reconstruct equivariant $K$-theory and that 
the fact that complex $K$-theory is obtained from a categorical construction 
is used in an essential way to define equivariant $K$-theory. Recently 
J. Lurie constructed not only equivariant versions  but also 
\emph{$2$-equivariant versions} of elliptic cohomology (see \cite[\S 5.4]{lu}).  This means that not 
only an action of a group can be incorporated in the definition of elliptic cohomology, 
but also an action of a \emph{$2$-group} (i.e. of a categorical group). Now, 
a $2$-group $G$ can not act in a very interesting way on an object in a $1$-category, as this
action would simply be induced by an action of the group $\pi_{0}(G)$. However, 
a $2$-group can definitely act in an interesting manner on an object in a $2$-category, 
since automorphisms of a given object naturally form a $2$-group. The
existence of $2$-equivariant version of elliptic cohomology therefore suggests again a
close relation between elliptic cohomology and $2$-categories. \\

\begin{center} \textbf{Towards a theory of categorical sheaves in algebraic geometry} \end{center}

The conclusion of the observation above is that there should exist an interesting
notion of \emph{categorical sheaves}, which are sheaves of categories
rather than sheaves of vector spaces, useful for a geometric description of objects underlying 
elliptic cohomology.
In this work we have been interested in this notion
independently of elliptic cohomology and in the context 
of algebraic geometry rather than topology. Although our final motivations
is to understand better elliptic cohomology we have found the theory of categorical 
sheaves in algebraic geometry interesting in its own and think that it deserves
a fully independent development. 

To be more precise, and to fix ideas, a categorical sheaf theory is required to satisfy the following conditions.

\begin{itemize}

\item To any scheme $X$ there exists a $2$-category $Cat(X)$, of categorical sheaves on $X$. 
The $2$-category $Cat(X)$ is expected to be a symmetric monoidal $2$-category. Moreover, 
we want $Cat(X)$ to be a categorification of the category $Mod(X)$ of sheaves of 
$\mathcal{O}_{X}$-modules on $X$, in the sense that there is a natural equivalence
between $Mod(X)$ and the category of endomorphisms of the unit object of $Cat(X)$.

\item The $2$-category $Cat(X)$ comes equiped with monoidal sub-$2$-categories
$Cat_{qcoh}(X)$, $Cat_{coh}(X)$ and $Cat_{parf}(X)$, which are categorifications of
the categories $QCoh(X)$, $Coh(X)$ and $Vect(X)$, of quasi-coherent sheaves, coherent sheaves, 
and vector bundles. The monoidal $2$-category $Cat_{parf}(X)$ is moreover expected to 
be rigid (i.e. every object is dualizable).

\item For a morphism $f : X \longrightarrow Y$ of schemes, there is a $2$-adjunction
$$f^{*} : Cat(Y) \longrightarrow Cat(X) \qquad Cat(Y) \longleftarrow Cat(X) : f_{*}.$$
The $2$-functors $f^{*}$ and $f_{*}$ are supposed to preserve the
sub-$2$-categories $Cat_{qcoh}(X)$, $Cat_{coh}(X)$ and $Cat_{parf}(X)$, under some finiteness conditions
on $f$. 

\item There exists a notion of short exact sequence in $Cat(X)$, which can be used in order
to define a Grothendieck group $K^{(2)}_{0}(X):=K_{0}(Cat_{parf}(X))$ (or more generally
a ring spectrum $K(Cat_{parf}(X))$). This Grothendieck group is called the \emph{secondary K-theory of $X$}
and is expected to possess the usual functorialities in $X$ (at least pull-backs and push-forwards along
proper and smooth morphisms). 

\item There exists a Chern character
$$K^{(2)}_{0}(X) \longrightarrow H^{(2)}(X),$$
for some \emph{secondary cohomology group} $H^{(2)}(X)$. This Chern character is expected to be
functorial for pull-backs and to satisfy some version of the Grothendieck-Riemann-Roch formula
for push-forwards.

\end{itemize}

As we will see in section \S 2, it is not clear how to develop a theory as above, and it seems
to us that a theory satisfying all the previous requirements cannot reasonably exist. One major observation 
here is that the situation becomes more easy to handle if the categories
$Mod(X)$, $QCoh(X)$, $Coh(X)$ and $Vect(X)$ are replaced by their derived analogs
$D(X)$, $D_{qcoh}(X)$, $D_{coh}^{b}(X)$ and $D_{parf}(X)$. Our wanted categorical sheaf theory 
should then rather a be \emph{derived categorical sheaf theory}, and is expected to 
satisfy the following conditions.

\begin{itemize}

\item To any scheme $X$ is associated a \emph{triangulated-$2$-category} $Dg(X)$, 
of \emph{derived categorical sheaves on $X$}. Here, by 
\emph{triangulated-$2$-category} we mean a $2$-category whose categories of morphisms
are endowed with triangulated structure in a way that the composition functors are bi-exacts.
The $2$-category $Dg(X)$ is expected to be a symmetric monoidal $2$-category, in way
which is compatible with the triangulated structure. Moreover, 
we want $Dg(X)$ to be a categorification of the derived category $D(X)$ of sheaves of 
$\mathcal{O}_{X}$-modules on $X$, in the sense that there is a natural triangulated equivalence
between $D(X)$ and the triangulated category of endomorphisms of the unit object of $Dg(X)$.

\item The $2$-category $Dg(X)$ comes equipped with monoidal sub-$2$-categories
$Dg_{qcoh}(X)$, $Dg_{coh}(X)$ and $Dg_{parf}(X)$, which are categorifications of
the derived  categories $D_{qcoh}(X)$, $D^{b}_{coh}(X)$ and $D_{parf}(X)$, of 
quasi-coherent complexes, bounded coherent sheaves, 
and perfect complexes. The monoidal $2$-category 
$Dg_{parf}(X)$ is moreover expected to be \emph{rigid} (i.e. every object 
is dualizable).

\item For a morphism $f : X \longrightarrow Y$ of schemes, there is a $2$-adjunction
$$f^{*} : Dg(Y) \longrightarrow Dg(X) \qquad Dg(Y) \longleftarrow Dg(X) : f_{*}.$$
The $2$-functors $f^{*}$ and $f_{*}$ are supposed to preserve the
sub-$2$-categories $Dg_{qcoh}(X)$, $Dg_{coh}(X)$ and $Dg_{parf}(X)$, under some finiteness conditions
on $f$. 

\item There exists a notion of short exact sequence in $Dg(X)$, which can be used in order
to define a Grothendieck group $K^{(2)}_{0}(X):=K_{0}(Dg_{parf}(X))$ (or more generally
a ring spectrum $K(Dg_{parf}(X))$). This Grothendieck group is called the \emph{secondary K-theory of $X$}
and is expected to possess the usual functorialities in $X$ (at least pull-backs and push-forwards along
proper and smooth morphisms). 

\item There exists a Chern character
$$K^{(2)}_{0}(X) \longrightarrow H^{(2)}(X),$$
for some \emph{secondary cohomology group} $H^{(2)}(X)$. This Chern character is expected to be
functorial for pull-backs and to satisfy some version of the Grothendieck-Riemann-Roch formula
for push-forwards.

\end{itemize}

The purpose of these notes is to give some ideas on how to define such 
triangulated-$2$-categories $Dg(X)$, how to define the secondary cohomology
$H^{(2)}(X)$ and how to define the Chern character. In order to do this, we will follow
closely one possible interpretation of the usual Chern character for vector bundles as being 
a kind of function on the loop space. \\

\begin{center} \textbf{The Chern character and the loop space} \end{center}

The Chern character we will construct for categorical sheaves is based on the
following interpretation of the usual Chern character. Assume that $X$ is
a smooth complex algebraic manifold (or more generally a complex algebraic stack) and 
that $V$ is a vector bundle on $X$. 
Let $\gamma : S^{1} \longrightarrow X$ be a \emph{loop} in $X$. We do not want to specify want 
we mean by a loop here, and the notions of loops we will use in the sequel is a rather unconventional
one (see \ref{dloop}). Whatever $\gamma$ truely is, we will think of it as a loop in $X$, at least intuitively.
We consider the pull-back $\gamma^{*}(V)$, which is a vector bundle on $S^{1}$. 
Because of the notion of loops we use this vector bundle is in fact locally constant 
on $S^{1}$, and thus is completely determined by a monodromy operator
$m_{\gamma}$ on the fiber $V_{\gamma(0)}$. The trace of $m_{\gamma}$ is a complex number, and
as $\gamma$ varies in $LX$ the loop space of $X$ (again the notion of loop space we use 
is unconventional) we obtain a function $Ch(V)$ on $LX$. This function can be seen 
to be $S^{1}$-equivariant and thus provides an element
$$Ch(V) \in \mathcal{O}(LX)^{S^{1}}.$$
Our claim is that, if the objects $S^{1}$ and $LX$ are defined correctly, then there is a
natural identification
$$\mathcal{O}(LX)^{S^{1}}\simeq H^{ev}_{DR}(X),$$
and that $Ch(V)$ is the usual Chern character with values in the algebraic de Rham cohomology of $X$.
The conclusion is that $Ch(V)$ can be seen as a $S^{1}$-equivariant function on $LX$. 

One enlightening example is when $X$ is $BG$ the quotient stack of a finite group $G$.
The our loop space $LBG$ is the quotient stack $[G/G]$, for the action of $G$ on itself by conjugation.
The space of functions on $LBG$ can therefore be identified with 
$\mathbb{C}(G)$, the space of class function on $G$. A vector bundle $V$ on $BG$ is nothing else
than a linear representation of $G$, and the function $Ch(V)$ constructed above is the class
function sending $g\in G$ to $Tr(g:V\rightarrow V)$. Therefore, the description of the Chern character above 
gives back the usual morphism $R(G) \longrightarrow \mathbb{C}(G)$ sending 
a linear representation to its class function.

Our construction of the Chern character for a categorical sheaf follows the same ideas. 
The interesting feature of the above interpretation of the Chern character is that 
it can be generalized to any setting for which traces of endomorphisms make sense. 
As we already mentioned, $Dg_{parf}(X)$ is expected to be a rigid monoidal $2$-category, and thus any
endomorphism of an object possesses a trace which is itself an object 
in $D_{parf}(X)\simeq End(1)$. Therefore, if we start with a categorical sheaf on $X$ and do the same construction 
as above, we get a sheaf (rather than a function) on $LX$, or more precisely an object
in $D_{parf}(LX)$. This sheaf is moreover invariant under the action of $S^{1}$ and 
therefore is an object in $D_{parf}^{S^{1}}(LX)$, the perfect $S^{1}$-equivariant derived category of $LX$.
This sheaf as itself a Chern character which is an element in $H^{S^{1}}_{DR}(LX)$, the $S^{1}$-equivariant
de Rham cohomology of $LX$. This element is by definition the Chern character of our categorical sheaf.
The Chern character should then expected to be a map
$$Ch : K^{(2)}_{0}(X) \longrightarrow H^{S^{1}}_{DR}(LX).$$

\begin{center} \textbf{Plan of the paper} \end{center}

The main purpose of this paper is to make precise all the terms of this construction. 
For this we will start by the definitions of the $2$-categories $Dg(X)$, $Dg_{qcoh}(X)$, 
 and $Dg_{parf}(X)$, but we do not try to define $Dg_{coh}(X)$ as the notion
of coherence in this categorical setting seems unclear at the moment. The objects in $Dg(X)$ will be 
certain sheaves of dg-categories on $X$ and our approach to the notion 
of categorical sheaves heavily relies on the homotopy theory of dg-categories 
recently studied in \cite{ta,to2}. In a second part we will recall briefly some ideas of derived algebraic
geometry and of derived schemes (and stacks) as introduced in \cite{hagii,lu2}. The loop space $LX$
of a scheme $X$ will then be defined as the derived mapping stack from $S^{1}=B\mathbb{Z}$
to $X$. We will argue that the ring of $S^{1}$-invariant functions on $LX$ can be naturally 
identified with $HC_{0}^{-}(X)$, the negative cyclic homology of $X$. We will also 
briefly explain how this can be used in order to interpret the Chern character with values in 
cyclic homology as we have sketched above. Finally, in a last 
part we will present the construction of our Chern character for categorical sheaves. 
One crucial point in this construction is to define an $S^{1}$-equivariant sheaf 
on the loop space $LX$. The construction of the sheaf itself is easy but 
the fact that it is $S^{1}$-equivariant  is a delicate question which we leave open 
in the present work (see \ref{conj1} and \ref{conj2}). Hopefully a detailed proof of the existence of this
$S^{1}$-equivariant sheaf will appear in a future work. \\

\section{Categorification of homological algebra and dg-categories}

In this section we present our triangulated-$2$-categories $Dg(X)$ of 
derived categorical sheaves on some scheme $X$. We will start by 
an overview of a rather standard way to categorify the theory of
modules over some base commutative ring using linear categories. As we will 
see the notion of $2$-vector spaces appear naturally in this setting as the dualizable
objects, exactly in the same way that the dualizable modules are the
projective modules of finite rank. After arguing that this notion of $2$-vector space is
too rigid a notion to allow for push-fowards, we will consider
dg-categories instead and show that they can be used in order to 
categorify homological algebra in a similar way as linear categories categorify
linear algebra. By analogy
with the case of modules and linear categories we will consider dualizable
objects as categorified versions of perfect complexes and notice that 
these are precisely the smooth and proper dg-categories studied in \cite{ks,tv}. 
We will finally define the $2$-categories $Dg(X)$, $Dg_{qcoh}(X)$
and $Dg_{parf}(X)$ for a general scheme $X$ by some gluing procedure. \\

Let $k$ be a commutative base ring. We let $Mod(k)$ be the category of 
$k$-modules, considered as a symmetric monoidal category for the tensor
product of modules. Recall that an object $M\in Mod(k)$ is said to be
dualizable if the natural morphism 
$$M\otimes M^{\vee} \longrightarrow \underline{Hom}(M,M)$$
is an isomorphism (here $\underline{Hom}$ denotes the $k$-module
of $k$-linear morphisms, and $M^{\vee}:=\underline{Hom}(M,k)$ is the dual module).
It is easy to see that $M$ is dualizable if and only if it is projective and of finite type over $k$. 

A rather standard way to categorify the category $Mod(k)$ is to consider
$k$-linear categories and Morita morphisms. We let $Cat(k)$ be the 
$2$-category whose objects are small $k$-linear categories. The category of
morphisms between two $k$-linear categories $A$ and $B$ in $Cat(k)$
is defined to be the category of all $A\otimes_{k}B^{op}$-modules (the composition
is obtained by the usual tensor product of bi-modules). The tensor 
product of linear categories endow $Cat(k)$ with a structure of a symmetric monoidal 
$2$-category for which $k$, the $k$-linear category freely generated by one object, 
is the unit. We have $End_{Cat(k)}(k)\simeq Mod(k)$, showing that 
$Cat(k)$ is a categorification of $Mod(k)$. To obtain a categorification 
of $Mod^{pft}(k)$, the category of projective $k$-modules of finite type, 
we consider the sub-$2$-category of $Cat(k)$ with the same objects but for which 
the category of morphisms from $A$ to $B$ is the full sub-category of 
the category of $A\otimes_{k}B^{op}$-modules whose objects are 
bi-modules $M$ such that for any $a\in A$ the $B^{op}$-module $M(a,-)$ is 
projective of finite type (i.e. a retract of a finit sum of representable
$B^{op}$-modules). We let $Cat^{c}(k) \subset Cat(k)$ be this sub-$2$-category, which 
is again a symmetric monoidal $2$-category for tensor product of linear categories. 
By definition we have $End_{Cat^{c}(k)}(k)\simeq Mod^{pft}(k)$. However, 
the tensor category $Mod^{pft}(k)$ is a rigid tensor category in the sense that 
every object is dualizable, but not every object in $Cat^{c}(k)$ is dualizable. We therefore
consider $Cat^{sat}(k)$ the full sub-$2$-category of dualizable objects in 
$Cat^{c}(k)$. Then, $Cat^{sat}(k)$ is a rigid monoidal $2$-category which is 
a categorification of $Mod^{pft}(k)$. It can be checked that 
a linear category $A$ is in $Cat^{sat}(k)$ if and only if it is equivalent
in $Cat(k)$ to an associative $k$-algebra $B$ (as usual considered as a
linear category with a unique object) satisfying the following two conditions.

\begin{enumerate}

\item The $k$-module $B$ is projective and of finite type over $k$.

\item For any associative $k$-algebra $A$, a $B\otimes_{k}A$-module $M$
is projective of finite type if and only if it is so as a $A$-module.

\end{enumerate}

These conditions are also equivalent to the following two conditions.

\begin{enumerate}

\item The $k$-module $B$ is projective and of finite type over $k$.

\item The $B\otimes_{k}B^{op}$-module $B$ is projective.

\end{enumerate}

When $k$ is a field, an object in $Cat^{sat}(k)$ 
is nothing else than a finite dimensional $k$-algebra $B$ which is
universally semi-simple (i.e. such that $B\otimes_{k}k'$ is semi-simple
for any field extension $k \rightarrow k'$). In general, 
an object in $Cat^{sat}(k)$ is a flat family of universally semi-simple finite dimensional
algebras over $Spec\, k$. In particular, if $k$ is an algebraically closed field
any object in $Cat^{sat}(k)$ is equivalent to $k^{n}$, or in other words is
a $2$-vector space of finite dimension in the sense of Kapranov-Voevodsky (see for instance
\cite{bdr}). For a general commutative
ring $k$, the $2$-category $Cat^{sat}(k)$ is a reasonable generalization of the
notion of $2$-vector spaces and can be called the \emph{$2$-category of $2$-vector bundles
on $Spec\, k$}. 

One major problem with this notion of $2$-vector bundles is the lack of 
push-forwards in general. For instance, let $X$ be a smooth and proper algebraic variety
over some algebraically closed field $k$. We can consider $\underline{Vect}$, the \emph{trivial $2$-vector bundle
of rank $1$ over $X$}, which is the stack in categories sending a Zariski open $U\subset X$
to the linear category $Vect(U)$ of vector bundles over $U$. The push-forward of this
trivial $2$-vector bundle along the structure morphism $X \longrightarrow Spec\, k$ 
is the $k$-linear category of global sections of $\underline{Vect}$, or in other words
the $k$-linear category $Vect(X)$ of vector bundles on $X$. This is
an object in $Cat(k)$, but is definitely not in $Cat^{sat}(k)$. The linear category 
$Vect(X)$ is big enough to convince anyone that it can not be \emph{finite dimensional}
in any reasonable sense. This shows that the global sections of 
a $2$-vector bundle on a smooth and proper variety is in general not 
a $2$-vector bundle over the base field, and that in general it is 
hopeless to expect a good theory of proper push-forwards in this setting. 

A major observation in this work is that considering a categorification of
$D(k)$ instead of $Mod(k)$, which is what we call a \emph{categorification
of homological algebra}, solves the problem mentioned above concerning push-forwards. 
Recall that a dg-category (over some base commutative ring $k$) is 
a category enriched over the category of complexes of $k$-modules (see \cite{ta}). For a dg-category $T$
we can define its category of $T$-dg-modules as well as as its derived category 
$D(T)$ by formally inverting quasi-isomorphisms between dg-modules (see \cite{to2}).
For two dg-categories $T_{1}$ and $T_{2}$ we can form their tensor product
$T_{1}\otimes_{k}T_{2}$, as well as their derived tensor product $T_{1}\otimes_{k}^{\mathbb{L}}T_{2}$
(see \cite{to2}). We now define a $2$-catgeory $Dg(k)$ whose objects are dg-categories and
whose category of morphisms from $T_{1}$ to $T_{2}$ is $D(T_{1}\otimes_{k}^{\mathbb{L}}T^{op}_{2})$. 
The composition of morphisms is defined using the derived tensor product
$$-\otimes_{T_{2}}^{\mathbb{L}} - : D(T_{1}\otimes^{\mathbb{L}}_{k}T_{2}^{op}) \times
D(T_{2}\otimes^{\mathbb{L}}_{k}T_{3}^{op}) \longrightarrow 
D(T_{1}\otimes^{\mathbb{L}}_{k}T_{3}^{op}).$$
Finally, the derived tensor product of dg-categories endows $Dg(k)$ with a structure of a 
symmetric monoidal $2$-category. 

The symmetric monoidal $2$-category $Dg(k)$ is a categorification of the derived category 
$D(k)$ as we have by definition
$$\underline{End}_{Dg(k)}(1)\simeq D(k).$$
To obtain a categorification of $D_{parf}(k)$, the perfect derived category, we
consider the sub-$2$-category $Dg^{c}(k)$ having the same of objects as
$Dg(k)$ itself but for which the category of morphisms from $T_{1}$ to $T_{2}$ 
in $Dg^{c}(k)$ is the full sub-category of $D(T_{1}\otimes_{k}^{\mathbb{L}}T_{2}^{op})$
of bi-dg-modules $F$ such that for all $t\in T_{1}$ the object
$F(t,-)\in D(T^{op})$ is compact (in the sense of triangulated categories, see \cite{ne}).
The symmetric monoidal structure on $Dg(k)$ restricts to a 
symmetric monoidal structure on $Dg^{c}(k)$, and we have
$$\underline{End}_{Dg^{c}(k)}(1)\simeq D_{parf}(k),$$
as an object of $D(k)$ is compact if and only if it is a perfect complex.
Finally, the symmetric monoidal $2$-category $Dg^{c}(k)$ is not rigid and we thus consider
$Dg^{sat}(k)$, the full sub-$2$-category consisting of rigid objects in $Dg^{c}(k)$. 
By construction, $Dg^{sat}(k)$ is a rigid symmetric monoidal $2$-category and we have
$$\underline{End}_{Dg^{sat}(k)}(1)\simeq D_{parf}(k).$$
The $2$-category $Dg^{sat}(k)$ will be our categorification of $D_{parf}(k)$ and its objects
should be thought as \emph{perfect derived categorical sheaves on the scheme $Spec\, k$}. 

It is possible to show that an object $T$ of $Dg^{c}(k)$ belongs to $Dg^{sat}(k)$ if 
and only if it is equivalent to an associative dg-algebra $B$, considered
as usual as a dg-category with a unique object, satisfying the following two conditions.
\begin{enumerate}
\item The underlying complex of $k$-modules of $B$ is perfect.
 
\item The object $B\in D(B\otimes_{k}^{\mathbb{L}}B^{op})$ is compact.
\end{enumerate}

In other words, a dg-category $T$ belongs to $Dg^{sat}(k)$ if and only
it is Morita equivalent to a smooth (condition $(2)$ above) and proper
(condition $(1)$ above) dg-algebra $B$. Such dg-categories are also often called
\emph{saturated} (see \cite{ks,tv}). 

As $Dg^{sat}(k)$ is a rigid symmetric monoidal $2$-category we can define, for any 
object $T$ a trace morphism
$$Tr : \underline{End}_{Dg^{sat}(k)}(T) \longrightarrow \underline{End}_{Dg^{sat}(k)}(1)\simeq D_{parf}(k).$$
Now, the category $\underline{End}_{Dg^{sat}(k)}(T)$ can be naturally identified with 
$D(T\otimes_{k}^{\mathbb{L}}T^{op})_{c}$, the full sub-category of $D(T\otimes_{k}^{\mathbb{L}}T^{op})$
of compact objects (here we use that $T$ is saturated and the results
of \cite[\S 2.2]{tv}). The trace morphism is then a functor
$$Tr : D(T\otimes_{k}^{\mathbb{L}}T^{op})_{c} \longrightarrow D_{parf}(k)$$
which can be seen to be isomorphic to the functor sending a bi-dg-module $M$ to its
Hochschild complex $HH(T,M) \in D_{parf}(k)$. In particular, 
the \emph{rank} of an object $T\in Dg^{sat}(k)$, which by definition is the trace of
its identity, is its Hochschild complex $HH(T)\in D_{parf}(k)$. \\

To finish this section we present the global versions of the $2$-categories $Dg(k)$ and
$Dg^{sat}(k)$ over some base scheme $X$. We let $ZarAff(X)$ be the small
site of affine Zariski open sub-schemes of $X$. We start to define 
a category $dg-cat(X)$ consisting of the following data
\begin{enumerate}
\item For any $Spec\, A = U \subset X$ in $ZarAff(X)$, 
a dg-category $T_{U}$ over $A$.
\item For any $Spec\, B = V \subset Spec\, A = U \subset X$ morphism 
in $ZarAff(X)$ a morphism of dg-categories over $A$
$$r_{U,V} : T_{U} \longrightarrow T_{V}.$$
\end{enumerate}
These data should moreover satisfy the equation $r_{V,W}\circ r_{U,V}=r_{U,W}$ for any
inclusion of affine opens $W \subset V \subset U\subset X$. The morphisms in 
$dg-cat(X)$ are defined in an obvious way as families of dg-functors commuting with 
the $r_{U,V}$'s. 

For 
$T\in dg-cat(X)$ we define a category $Mod(T)$ of $T$-dg-modules in the following way.
Its objects consist of the following data
\begin{enumerate}
\item For any $Spec\, A = U \subset X$ in $ZarAff(X)$, 
a $T_{U}$-dg-module $M_{U}$.
\item For any $Spec\, B = V \subset Spec\, A = U \subset X$ morphism 
in $ZarAff(X)$ a morphism of $T_{U}$-dg-modules 
$$m_{U,V} : M_{U} \longrightarrow r_{U,V}^{*}(M_{V}).$$
\end{enumerate}
These data should moreover satisfy the usual cocyle equation for 
$r_{U,V}^{*}(m_{V,W})\circ m_{U,V}=m_{U,W}$. Morphisms in 
$Mod(T)$ are simply defined as families of morphisms of dg-modules commuting with the
$m_{U,V}$'s. Such a morphism $f : M \longrightarrow M'$ in $Mod(T)$ is 
a quasi-isomorphism if it is a stalkwise quasi-isomorphism (note that 
$M$ and $M'$ are complexes of presheaves of $\mathcal{O}_{X}$-modules). We denote
by $D(T)$ the category obtained from $Mod(T)$ by formally inverting these quasi-isomorphisms. 

We now define a $2$-category$Dg(X)$  whose objects are the objects of $dg-cat(X)$, and whose
category of morphisms from $T_{1}$ to $T_{2}$ is $D(T_{1}\otimes_{\mathcal{O}_{X}}^{\mathbb{L}}T_{2}^{op})$
(we pass on the technical point of defining this derived tensor product over $\mathcal{O}_{X}$, one
possibility being to endow $dg-cat(X)$ with a model category structure and to use a cofibrant replacement).
The compositions of morphisms in $Dg(X)$ is given by the usual derived tensor product. The derived tensor
product endows $Dg(X)$ with a structure of a symmetric monoidal $2$-category and we have by construction
$$\underline{End}_{Dg(X)}(1)\simeq D(X),$$
where $D(X)$ is the (unbounded) derived category of all $\mathcal{O}_{X}$-modules on $X$. 

\begin{df}\label{d1}
\begin{enumerate}
\item 
\emph{An object $T\in Dg(X)$ is} quasi-coherent \emph{if for any inclusion of affine open subschemes
$$Spec\, B = V \subset Spec\, A = U \subset X$$
the induced morphism
$$T_{U}\otimes_{A}B \longrightarrow T_{V}$$
is a Morita equivalence of dg-categories.} 
\item \emph{A morphism $T_{1} \longrightarrow T_{2}$, corresponding to an object 
$M\in D(T_{1}\otimes_{\mathcal{O}_{X}}^{\mathbb{L}}T_{2}^{op})$, is called} quasi-coherent
\emph{if its underlying complex of $\mathcal{O}_{X}$-modules is with quasi-coherent cohomology sheaves}.\\
\end{enumerate}

\emph{The sub-$2$-category of $Dg(X)$ consisting of quasi-coherent objects and quasi-coherent 
morphisms is denoted by $Dg_{qcoh}(X)$. It is called the }
$2$-category of quasi-coherent derived categorical sheaves on $X$.
\end{df}

Let $T_{1}$ and $T_{2}$ be two objects in $Dg_{qcoh}(X)$. We consider the full sub-category of
$D(T_{1}\otimes_{\mathcal{O}_{X}}^{\mathbb{L}}T_{2}^{op})$ consisting of 
objects $M$ such that for any Zariski open $Spec\, A = U \subset X$
and any object $x\in (T_{1})_{U}$, the induced dg-module
$M(a,-)\in D((T_{2}^{op})_{U})$
is compact. This defines a sub-$2$-category of 
$Dg_{qcoh}(X)$, denoted by $Dg_{qcoh}^{c}(X)$ and will be called
the sub-$2$-category of \emph{compact morphisms}. The symmetric monoidal structure 
on $Dg(X)$ restricts to a symmetric monoidal structure on $Dg_{qcoh}(X)$ and on 
$Dg_{qcoh}^{c}(X)$.

\begin{df}\label{d2}
\emph{The} $2$-category of perfect derived categorical sheaves \emph{is the full sub-$2$-catgeory 
of $Dg_{qcoh}^{c}(X)$ consisting of dualizable objects. It is denoted by 
$Dg_{parf}(X)$. }
\end{df}

By construction, $Dg_{parf}(X)$ is a symmetric monoidal $2$-category with 
$$\underline{End}_{Dg_{parf}(X)}(1)\simeq D_{parf}(X).$$

It is possible to show that an object $T\in D_{qcoh}(X)$ belongs to 
$Dg_{parf}(X)$ if and only if for any affine Zariski open subscheme
$Spec\, A = U \subset X$, the dg-category $T_{U}$ is saturated (i.e.
belongs to $Dg^{sat}(A)$). 

For a morphism of schemes $f : X \longrightarrow Y$ it is possible to define
a $2$-adjunction
$$f^{*} : Dg(Y) \longrightarrow Dg(X) \qquad Dg(Y) \longleftarrow Dg(X) : f_{*}.$$
Moreover, $f^{*}$ preserves quasi-coherent objects, quasi-coherent morphisms, 
as well as the sub-$2$-categories of compact morphisms and perfect objects. When the morphism
$f$ is quasi-compact and quasi-separated, we think that it is possible to prove 
that $f_{*}$ preserves quasi-coherent objects and quasi-coherent morphisms as well as
the sub-$2$-category of compact morphisms. We also guess that 
$f_{*}$ will preserve perfect objects when $f$ is smooth and proper, but this would
require a precise investigation. As a typical example, the direct image of the unit 
$1\in Dg_{parf}(X)$ by $f$ is the presheaf of dg-categories sending $Spec\, A = U \subset Y$ to 
the dg-category $L_{parf}(X\times_{Y}U)$ of perfect complexes over the scheme
$f^{-1}(U)\simeq X\times_{Y}U$.  When $f$ is smooth and proper it is known that 
the dg-category $L_{parf}(X\times_{Y}U)$ is in fact saturated (see \cite[\S 8.3]{to2}). This shows that 
$f_{*}(1)$ is a perfect derived sheaf on $Y$ and provides an evidence that 
$f_{*}$ preserves perfect object. These functoriality statements will be considered in more details
in a future work. 

\section{Loop spaces in derived algebraic geometry}

In this section we present a version of the loop space of a scheme (or more generally 
of an algebraic stack) based on derived algebraic geometry. For us the circle
$S^{1}$ is defined to be the quotient stack $B\mathbb{Z}$, where 
$\mathbb{Z}$ is considered as a constant sheaf of groups. For any scheme $X$, 
the mapping stack $\mathbf{Map}(S^{1},X)$ is then equivalent to $X$, as the coarse moduli 
space of $S^{1}$ is simply a point. In other words, with this definition of the circle there are
no interesting loops on a scheme $X$. However, we will explain in the sequel that 
there exists an interesting \emph{derived mapping stack} $\mathbb{R}\mathbf{Map}(S^{1},X)$, 
which is now a derived scheme and which is non-trivial. This derived mapping stack
will be our loop space. In this section we recall briefly the
notions from derived algebraic geometry needed in order to define 
the object $\mathbb{R}\mathbf{Map}(S^{1},X)$. We will also explain the relation between 
the cohomology of $\mathbb{R}\mathbf{Map}(S^{1},X)$ and cyclic homology of the scheme $X$. \\

Let $k$ be a base commutative ring and denote by $\mathbf{Sch}_{k}$, resp. $\mathbf{St}_{k}$,
the category of schemes over $k$ and the model category of stacks over $k$ (\cite[2.1.1]{hagii} or \cite[\S 2, 3]{seattle}),
for the \'etale topology. We recall that the homotopy category $\mathrm{Ho}(\mathbf{St}_{k})$
contains as full sub-categories the category of sheaves of sets on $\mathbf{Sch}_{k}$ as well as the
$1$-truncation of the $2$-category of stacks in groupoids (in the sense of \cite{lm}). In particular, 
the homotopy category of stacks $\mathrm{Ho}(\mathbf{St}_{k})$ contains the category of
schemes and of Artin stacks as full sub-categories. In what follows we will always consider these two categories 
as embedded in $\mathrm{Ho}(\mathbf{St}_{k})$. Finally, recall that the category 
$\mathrm{Ho}(\mathbf{St}_{k})$ possesses internal Hom's, that will be denoted by 
$\mathbf{Map}$. 

As explained in \cite[Ch. 2.2]{hagii} (see also \cite[\S 4]{seattle} for an overview), 
there is also a model category $\mathbf{dSt}_{k}$ of derived stacks
over $k$ for the strong \'etale topology. The derived affine objects are simplicial $k$-algebras, and the model category
of simplicial $k$-algebras will be denoted by $\mathbf{salg}_{k}$. 
The opposite (model) category is denoted by $\mathbf{dAff}_{k}$. Derived stacks can be identified with objects in the homotopy category
$\mathrm{Ho}(\mathbf{dSt}_{k})$ which in turn can be identified with the full subcategory of the homotopy category of
simplicial presheaves on $\mathbf{dAff}_{k}$ whose objects are weak equivalences' preserving simplicial presheaves $F$
having strong \'etale descent i.e.
such that, for any \'etale homotopy hypercover $U_{\bullet}\rightarrow X$ in $\mathbf{dAff}_{k}$ (\cite[Def. 3.2.3, 4.4.1]{hagi}), the canonical map 
$$F(X)\longrightarrow \mathrm{holim}F(U_{\bullet})$$ is a weak equivalence of simplicial sets. The derived Yoneda functor
induces a fully faithful functor on homotopy categories $$\mathbb{R}\mathrm{Spec}:\mathrm{Ho}(\mathbf{dAff}_{k})\hookrightarrow \mathrm{Ho}(\mathbf{dSt}_{k})\; : A\mapsto (\,\mathbb{R}\mathrm{Spec}(A)\; : B\mapsto \mathrm{Map}_{\mathbf{salg}_{k}}(A,B) \,),$$
where $\mathrm{Map}_{\mathbf{salg}_{k}}$ denotes the mapping spaces of the model 
category $\mathbf{salg}_{k}$ (therefore $\mathrm{Map}_{\mathbf{salg}_{k}}(A,B)\simeq \underline{Hom}(Q(A),B)$, 
where $\underline{Hom}$ denotes the natural simplicial Hom's of $\mathbf{salg}_{k}$
and $Q(A)$ is a cofibrant model for $A$).
Those derived stacks belonging to the essential image of $\mathbb{R}\mathrm{Spec}$ will be called \textit{affine} derived stacks. 

The category $\mathrm{Ho}(\mathbf{dSt}_{k})$ of derived stacks has a lot of important properties. First of all, being the homotopy category of a model category, it has derived colimits and limits (denoted as $\mathrm{hocolim}$ and $\mathrm{holim}$). In particular, given any pair of maps $F\rightarrow S$ and $G\rightarrow S$ between derived stacks, there is a derived fiber product stack $\mathrm{holim}(F\rightarrow S \leftarrow G)\equiv F\times_{S} ^{h} G$. As our base ring $k$ is not assumed to be a field, 
the direct product in the model category $\mathbf{dSt}_{k}$ is not exact and should also be derived. 
The derived direct product of two derived stacks $F$ and $G$ will be denoted by $F\times^{ h}G$. This derived product 
is the categorical product in the homotopy category  $\mathrm{Ho}(\mathbf{dSt}_{k})$.
The category  $\mathrm{Ho}(\mathbf{dSt}_{k})$ also admits internal Hom's, i.e. for any pair of derived stacks $F$ and $G$ there is a derived mapping stack denoted as $$\mathbb{R}\mathbf{Map}(F,G)$$ with the property that 
$$[F,\mathbb{R}\mathbf{Map}(G,H)]\simeq 
[F\times^{h} G,H],$$ 
functorially in $F$, $G$ and $H$.

The inclusion functor $j$ of commutative $k$-algebras into $\mathbf{salg}_{k}$ (as constant simplicial algebras) induces
a pair $(i,t_0)$ of (left,right) adjoint functors $$t_{0}:=j^* \,: \mathrm{Ho}(\mathbf{dSt}_{k})\rightarrow \mathrm{Ho}(\mathbf{St}_{k}) \qquad i:=\mathbb{L}j_!: \mathrm{Ho}(\mathbf{St}_{k})\rightarrow \mathrm{Ho}(\mathbf{dSt}_{k}).$$
It can be proved that $i$ is fully faithful. In particular we can, and will, view any stack as a derived stack (we will
most of the time omit to mention the functor $i$ and consider $\mathrm{Ho}(\mathbf{St}_{k})$ as embedded
in $\mathrm{Ho}(\mathbf{dSt}_{k})$). 
The truncation functor $t_0$ acts on affine derived stacks as $t_0 (\mathbb{R}\mathrm{Spec}(A))=\mathrm{Spec}(\pi_0 A)$. 
It is important to note that the inclusion functor $i$ does not preserve derived internal hom's nor derived fibered products. 
This is a crucial point in derived algebraic geometry: derived tangent spaces and derived fiber products of usual schemes or stacks are really derived objects. The derived tangent space of an Artin stack viewed as a
 derived stack via $i$ is the dual of its cotangent complex while the derived fiber product of, say, two affine schemes viewed as two derived stacks is given by the derived tensor product of the corresponding commutative algebras $$i(\mathrm{Spec}\,S) \times_{i(\mathrm{Spec}\,R)} ^{h}i(\mathrm{Spec}\,T) \simeq \mathbb{R}\mathrm{Spec}\,(S\otimes_{R}^{\mathbb{L}}T).$$\\

Both for stacks and derived stacks there is a notion of being \textit{geometric} (\cite[1.3, 2.2.3]{hagii}), depending, among other things, on the choice of a notion of smooth morphism between the affine pieces. For morphisms of commutative $k$-algebras this is the usual notion of smooth morphism, while in the derived case, a morphism $A\rightarrow B$ of simplicial $k$-algebras is said to be strongly smooth if the induced map $\pi_0 A \rightarrow\pi_0 B$ is a smooth morphism of commutative rings, and $\pi_{*}A\otimes_{\pi_{0}A} \pi_0 B\simeq \pi_{*}B$. The notion of geometric stack is strictly related to the notion of Artin stack (\cite[Prop. 2.1.2.1]{hagii}). Any geometric derived stack has a cotangent complex (\cite[Cor. 2.2.3.3]{hagii}). Moreover, both functors $t_0$ and $i$ preserve geometricity.\\

Let $\mathrm{B}\mathbb{Z}$ be the classifying stack of the constant group scheme $\mathbb{Z}$. We view $\mathrm{B}\mathbb{Z}$
as an object of $\mathbf{St}_{k}$, i.e. as the stack associated to the constant simplicial presheaf $$\mathrm{B}\mathbb{Z}: \mathbf{alg}_{k}\rightarrow \mathbf{SSets} \, , R\mapsto \mathrm{B}\mathbb{Z} ,$$ where by abuse of notations, we have also denoted as $\mathrm{B}\mathbb{Z}$ the classifying simplicial set, i.e. the nerve of the (discrete) group $\mathbb{Z}$. Such a nerve is naturally a pointed simplicial set, and we call $0$  that point.

\begin{df}\label{dloop}\emph{Let $X$ be a derived stack over $k$. The} derived loop stack \emph{of $X$ is the derived stack $$\mathrm{L}X:= \mathbb{R}\mathbf{Map}(\mathrm{B}\mathbb{Z}, X).$$}
\end{df} 

We will be mostly interested in the case where $X$ is a scheme or an algebraic (underived) stack.
Taking into account the homotopy equivalence 
$$\mathrm{B}\mathbb{Z} \simeq S^{1}\simeq *\coprod^{h}_{*\coprod *} * \, ,$$ 
we see that we  have 
$$\mathrm{L}X\simeq X\times_{X\times^{h} X}^{h}X,$$ 
where $X$ maps to $X\times^{h} X$ diagonally and the homotopy fiber product is taken in $\mathbf{dSt}_{k}$.
Evaluation at $0\in \mathrm{B}\mathbb{Z}$ yields a canonical map of derived stacks 
$$p:\mathrm{L}X \longrightarrow X.$$ 
On the other hand, since the limit maps canonically to the homotopy limit, we get a canonical morphism of derived stacks $X \rightarrow \mathrm{L}X$, a section of $p$, describing $X$ as the ``constant loops'' in $\mathrm{L}X$. 

If $X$ is an affine scheme over $k$,  $X=\mathrm{Spec}\, A$ with $A$ a commutative $k$-algebra, we get that 
$$\mathrm{L}X\simeq \mathbb{R}\mathrm{Spec}(A\otimes_{A\otimes^{\mathbb{L}} A}^{\mathbb{L}}A),$$ 
where the derived tensor product is taken in the model category $\mathbf{salg}_{k}$. One way to rephrase this is by saying that ``functions'' on $\mathrm{L}\mathrm{Spec}\, A$ \textit{are} Hochschild homology classes of $A$ with values in $A$ itself. 
Precisely, we have 
$$\mathcal{O}(\mathrm{L}X):=\mathbb{R}\underline{Hom}(\mathrm{L}X,\mathbb{A}^{1})
\simeq \mathrm{HH}(A,A),$$ 
where $\mathrm{HH}(A,A)$ is the simplicial set obtained from the complex of Hochschild homology of $A$
by the Dold-Kan correspondence, and $\mathbb{R}\underline{Hom}$ denotes the
natural  enrichment of $\mathrm{Ho}(\mathbf{dSt}_{k})$
into $\mathrm{Ho}(SSet)$. When $X$ is a general scheme then 
$\mathcal{O}(\mathrm{L}X)$ can be identified with the Hochschild homology complex of $X$, and we have
$$\pi_{i}(\mathcal{O}(\mathrm{L}X))\simeq HH_{i}(X).$$
In particular, when $X$ is a smooth and $k$ is of characteristic zero, the Hochschild-Kostant-Rosenberg theorem 
implies that 
$$\pi_{0}(\mathcal{O}(\mathrm{L}X))\simeq \oplus_{i} H^{i}(X,\Omega_{X/k}^{i}).$$

The stack $S^{1}=\mathrm{B}\mathbb{Z}$ is a group stack, and it acts naturally on 
$\mathrm{L}X$ for any derived stack $X$ by ``rotating the loops''. More precisely, 
there is a model category $\dst^{S^{1}}$, or $S^{1}$-equivariant stacks, and 
$\mathrm{L}X$ is naturally an object in the homotopy category $\mathrm{Ho}(\dst^{S^{1}})$. 
This way, the simplicial algebra of functions $\mathcal{O}(\mathrm{L}X)$ is naturally 
an $S^{1}$-equivariant simplicial algebra, and thus can also be considered
as an $S^{1}$-equivariant complex or in other words as an object 
in $D^{S^{1}}(k)$, the $S^{1}$-equivariant derived category of $k$. The category $D^{S^{1}}(k)$ is
also naturally equivalent to $D(k[\epsilon])$, the derived category of the dg-algebra $k[\epsilon]$
freely generated by an element $\epsilon$ of degree $-1$ and with $\epsilon^{2}=0$. The derived category
$D^{S^{1}}(k)$ is thus naturally equivalent to the derived category of mixed complexes (see
\cite{lo}) (multiplication by $\epsilon$ providing the second differential). When $X=\mathrm{Spec}\, A$ is an affine scheme,  
$\mathcal{O}(\mathrm{L}X)$ can be identified, as a mixed complex, with 
the Hoschschild complex $\mathrm{HH}(A)$, with its canonical mixed complex structure. As a consequence, we have
$$\pi_{i}(\mathcal{O}(\mathrm{L}X)^{hS^{1}})\simeq HC_{i}^{-}(A),$$
where $HC^{-}$ denotes negative cyclic homology and $K^{hS^{1}}$ denotes the simplicial set of
homotopy fixed points of an $S^{1}$-equivariant simplicial set $K$. In other words, 
there is a natural identification between $S^{1}$-invariant functions on $\mathrm{L}X$ and
negative cyclic homology of $X$. This statement of course can be generalized  to the case of
a scheme $X$.

\begin{prop}\label{p1}
For a scheme $X$ we have
$$\pi_{0}(\mathcal{O}(\mathrm{L}X)^{hS^{1}})\simeq HC_{0}^{-}(X),$$
where the right hand side denotes negative cyclic homology of the scheme $X$.
\end{prop}

When $k$ is of characteristic zero and $X$ is smooth over $k$ then proposition \ref{p1} 
states that 
$$\pi_{0}(\mathcal{O}(\mathrm{L}X)^{hS^{1}})\simeq H_{DR}^{ev}(X/k).$$
The even part of de Rham cohomology of $X$ can be identified with 
$S^{1}$-equivariant functions on the derived loop space $\mathrm{L}X$. This fact can also be
generalized to the case where $X$ is a smooth Deligne-Mumford stack
over $k$ (again assumed to be of characteristic zero), but the right hand side
should rather be replaced by the (even part of) de Rham orbifold cohomology of $X$, which is
the de Rham cohomology of the inertia stack $IX\simeq t_{0}(\mathrm{L}X)$. 

To finish this part, we would like to mention that the construction of the Chern character
for vector bundles
we suggested in section \S 1 can now be made precise, and through the identification 
of proposition \ref{p1} this Chern character coincides with the usual one. We start with
a vector bundle $V$ on $X$ and we consider its pull-back
$p^{*}(V)$ on $\mathrm{L}X$, which is a vector bundle on the derived scheme 
$\mathrm{L}X$. This vector bundle $p^{*}(V)$ comes naturally equipped with 
an automorphism $u$. This follows by considering the evaluation
morphism $\pi : S^{1}\times \mathrm{L}X \longrightarrow X$, and the
vector bundle $\pi^{*}(V)$. As $S^{1}=B\mathbb{Z}$, a vector bundle on 
$S^{1}\times \mathrm{L}X$ consists precisely of a vector bundle on 
$\mathrm{L}X$ together with an action of $\mathbb{Z}$, or in other words together with 
an automorphism. We can then consider the trace of $u$, which is an 
element in $\pi_{0}(\mathcal{O}(\mathrm{L}X))\simeq HH_ {0}(X)$. A difficult issue here 
is to argue that this function $Tr(u)$ has a natural refinement to an $S^{1}$-invariant function
$Tr(u) \in \pi_{0}(\mathcal{O}(\mathrm{L}X)^{h S^{1}})\simeq HC^{-}_ {0}(X)$, which is the
Chern character of $V$. The $S^{1}$-invariance of $Tr(u)$ will be studied in 
a future work, and we refer to our last section below, for some comments about how this 
would follow from the general theory of rigid tensor $\infty$-categories.

\section{Construction of the Chern character}

We are now ready to sketch the construction of our Chern character for a derived categorical sheaf. 
This construction simply follows the lines we have just sketched for vector bundles. We will meet the same 
difficult issue of the existence of an $S^{1}$-invariant refinement of the trace, and we will leave this
question as an conjecture. However, in the next section we will explain how this would follow from a very general fact 
about rigid monoidal $\infty$-categories. \\

Let $T \in Dg_{parf}(X)$ be a perfect derived categorical sheaf on some scheme $X$ (or more generally
on some algberaic stack $X$). We consider the natural morphism $p : \mathrm{L}X \longrightarrow X$ and
we consider $p^{*}(T)$, which is a perfect derived categorical sheaf on $\mathrm{L}X$. We have not defined
the notions of categorical sheaves on derived schemes or derived stacks but this is rather straightforward. 
As in the case of vector bundles explained in the last section, 
the object $p^{*}(T)$ comes naturally equipped with an autoequivalence $u$. This again follows from the fact that 
a derived categorical sheaf on $S^{1}\times \mathrm{L}X$ is the same thing as 
a derived categorical sheaf on $\mathrm{L}X$ together with an autoequivalence. We consider the trace of 
$u$ in order to get a perfect complex on the derived loop space
$$Tr(u)\in \underline{End}_{Dg_{parf}(\mathrm{L}X)}(1)=D_{parf}(\mathrm{L}X).$$
The main technical difficulty here is to show that $Tr(u)$ possesses a natural 
lift as an $S^{1}$-equivariant complex on $\mathrm{L}X$. We leave this as a conjecture.

\begin{conj}\label{conj1}
The complex $Tr(u)$ has a natural lift 
$$Tr^{S^{1}}(u)\in D^{S^{1}}_{parf}(\mathrm{L}X),$$ 
where $D^{S^{1}}_{parf}(\mathrm{L}X)$ is the $S^{1}$-equivariant perfect derived category of
$\mathrm{L}X$. 
\end{conj}

The above conjecture is not very precise as the claim is not that 
a lift simply exists, but rather than there exists a natural one. One of the difficulty in the conjecture above is that 
it seems difficult to characterize the required lift by some specific properties. We will see however that 
the conjecture can be reduced to a general conjecture about rigid monoidal $\infty$-categories. 

Assuming conjecture \ref{conj1}, we have $Tr^{S^{1}}(u)$ and we now consider its
class in the Grothendieck group of the triangulated category $D^{S^{1}}_{parf}(\mathrm{L}X)$. This is
our definition of the categorical Chern character of $T$.

\begin{df}\label{d4}
\emph{The} categorical Chern character of $T$ \emph{is}
$$Ch^{cat}(T):=[Tr^{S^{1}}(u)] \in K^{S^{1}}_{0}(\mathrm{L}X):=K_{0}(D^{S^{1}}_{parf}(\mathrm{L}X)).$$
\end{df}

The categorical Chern character $Ch^{cat}(T)$ can be itself refined into a 
\emph{cohomological Chern character} by using now the $S^{1}$-equivariant Chern character 
map for $S ^{1}$-equvariant perfect complexes on $\mathrm{L}X$. We skip some technical details here but the final result
is an element
$$Ch^{coh}(T):=Ch^{S^{1}}(Ch^{cat}(T)) \in \pi_{0}(\mathcal{O}(\mathrm{L}^{(2)}X)^{h (S^{1}\times S^{1})}),$$
where $\mathrm{L}^{(2)}X:=\mathbb{R}\mathbf{Map}(S^{1}\times S^{1},X)$ is now the
derived double loop space of $X$. The space $\pi_{0}(\mathcal{O}(\mathrm{L}^{(2)}X)^{h (S^{1}\times S^{1})})$
can reasonably be called the \emph{secondary negative cyclic homology of $X$} and 
should be thought (and actually is) the $S^{1}$-equivariant negative cyclic homology of
$\mathrm{L}X$. We therefore have
$$Ch^{coh}(T) \in HC_{0}^{-,S^{1}}(\mathrm{L}X).$$

\begin{df}\label{d5}
\emph{The} cohomological Chern character of $T$ \emph{is} 
$$Ch^{coh}(T):=Ch^{S^{1}}(Ch^{cat}(T)) \in HC_{0}^{-,S^{1}}(\mathrm{L}X):=\pi_{0}(\mathcal{O}(\mathrm{L}^{(2)}X)^{h S^{1}\times S^{1}})$$
defined above.
\end{df}

Obviously, it is furthermore expected that the constructions $T\mapsto Ch^{cat}(T)$ and $T \mapsto Ch^{coh}(T)$
satisfy standard properties such as additivity,
multiplicativity and functoriality with respect to pull-backs. The most general version of our Chern character map should be a morphism 
of commutative ring spectra
$$Ch^{cat} : Kg(X) \longrightarrow K^{S^{1}}(\mathrm{L}X),$$
where $Kg(X)$ is a ring spectrum constructed using a certain Waldhausen 
category of perfect derived categorical sheaves on $X$ and 
$K^{S^{1}}(\mathrm{L}X)$ is the $K$-theory spectrum of $S^{1}$-equivariant perfect complexes
on $\mathrm{L}X$. This aspect of the Chern character will again
be investigated in more details in a future work.

\section{Final comments}

\begin{center} \textbf{On $S^{1}$-equivariant trace maps} \end{center}

Our conjecture \ref{conj1} can be clarified using the language of higher categories. 
Recall that a $(1,\infty)$-category is an $\infty$-category in which all 
$n$-morphisms are invertible (up to higher morphisms) as soon as $n>1$. There exist
several well behaved models for the theory of $(1,\infty)$-categories, such as
simplicially enriched categories, quasi-categories, Segal categories and Rezk's spaces. We refer to 
\cite{be} for an overview of these various notions. What we will say below can be done in any of these theories, 
but, to fix ideas, we will work with $\mathbb{S}$-categories (i.e. simplicially enriched categories). 

We will be using $Ho(\mathbb{S}-Cat)$ the homotopy category of $\mathbb{S}$-categories, which is
the category obtained from the category of $\mathbb{S}$-categories and $\mathbb{S}$-functors by inverting
the (Dwyer-Kan) equivalences (a mixture between weak equivalences of simplicial sets and categorical
equivalences). An important property of $Ho(\mathbb{S}-Cat)$  is that it is cartesian closed (see \cite{to2}
for the corresponding statement for dg-categories whose proof is similar). 
In particular, for two $\mathbb{S}$-categories $T$ and $T'$ we can construct 
an $\mathbb{S}$-category $\mathbb{R}\underline{Hom}(T,T')$ with the property that 
$$[T'',\mathbb{R}\underline{Hom}(T,T')] \simeq [T''\times T,T'],$$
where $[-,-]$ denote the Hom sets of $Ho(\mathbb{S}-Cat)$.  Any 
$\mathbb{S}$-category $T$ gives rise to a genuine category $[T]$ with the same objects
and whose sets of morphisms are the connected components of the simplicial sets of morphisms
of $T$.

We let $\Gamma$ be the category of pointed finite sets and pointed maps. 
The finite set $\{0,\dots,n\}$ pointed at $0$ will be denoted by $n^{+}$.
Now, a \emph{symmetric monoidal
$\mathbb{S}$-category} $T$ is a functor
$$T : \Gamma \longrightarrow \mathbb{S}-Cat$$
such that for any $n\geq 0$ the so-called \emph{Segal morphism}
$$T(n^{+}) \longrightarrow T(1^{+})^{n},$$
induced by the various projections $n^{+} \rightarrow 1^{+}$ sending $i\in \{1,\dots,n\}$ to $1$ and everything else
to $0$, is an equivalence of $\mathbb{S}$-categories. The full sub-category of the homotopy category 
of functors $Ho(\mathbb{S}-Cat^{\Gamma})$ consisting of symmetric monoidal 
$\mathbb{S}$-categories will be denoted by $Ho(\mathbb{S}-Cat^{\otimes})$.  As
the category $Ho(\mathbb{S}-Cat)$ is a model for the homotopy category of
$(1,\infty)$-categories, the category $Ho(\mathbb{S}-Cat^{\otimes})$ is a model for 
the homotopy category of symmetric monoidal $(1,\infty)$-categories. For 
$T\in Ho(\mathbb{S}-Cat^{\otimes})$ we will again use $T$ to denote its underlying 
$\mathbb{S}$-category $T(1^{+})$. The $\mathbb{S}$-category $T(1^{+})$ 
has a natural structure of a commutative monoid in $Ho(\mathbb{S}-Cat)$. This
monoid structure will be denoted by $\otimes$. 

We say that a symmetric monoidal $\mathbb{S}$-category $T$ is \emph{rigid} if 
for any object $x\in T$ there is an objetc $x^{\vee} \in T$ and a morphism
$1 \rightarrow x\otimes x^{\vee}$ such that for any pair of 
objects $y,z\in T$, the induced morphism of simplicial sets
$$T(y\otimes x,z) \longrightarrow T(y\otimes x \otimes x^{\vee},z\otimes x^{\vee}) \longrightarrow T(y,z\otimes x^{\vee})$$
is an equivalence. In particular, the identity of $x^{\vee}$ provides a trace morphism 
$x\otimes x^{\vee} \rightarrow 1$ ($y=x^{\vee}$, $z=1$). Therefore, 
for any rigid symmetric monoidal $\mathbb{S}$-category $T$ and an object $x\in T$ we can define
a trace morphism
$$Tr_{x} : T(x,x) \simeq T(1,x\otimes x^{\vee}) \longrightarrow T(1,1).$$
Let $T$ be a fixed rigid symmetric monoidal $\mathbb{S}$-category and 
$S^{1}=B\mathbb{Z}$ be the groupoid with a unique object with
$\mathbb{Z}$ as automorphism group. The category $S^{1}$ is an abelian group object
in categories and therefore can be considered as a group object in $\mathbb{S}$-categories.
The $\mathbb{S}$-category of functors $\mathbb{R}\underline{Hom}(S^{1},T)$
is denoted by $T(S^{1})$, and is equiped with a natural action 
of $S^{1}$. We consider the sub-$\mathbb{S}$-category of 
invertible (up to homotopy) morphisms in $T(S^{1})$ whose classifying space
is an $S^{1}$-equivariant simplicial set. We denote this simplicial set 
by $LT$ (``$L$'' stands for ``loops''). It is possible to put all the trace morphisms $Tr_{x}$ defined above
into a morphism of simplicial sets (well defined in $Ho(SSet)$)
$$Tr : LT \longrightarrow T(1,1).$$
Note that the connected components of $LT$ are in one to one correspondence 
with the set of equivalences classes of pairs $(x,u)$, consisting of an object
$x$ in $T$ and an autoequivalence $u$ of $x$. The morphism $Tr$ is such that 
$Tr(x,u)=Tr_{x}(u)\in \pi_{0}(T(1,1))$. We are in fact convinced that 
the trace map $Tr$ can be made equivariant for the action of $S^{1}$ on $LT$, functorially in 
$T$. To make a precise conjecture, we consider $\mathbb{S}-Cat^{rig}$ the category 
of all rigid symmetric monoidal $\mathbb{S}$-categories (note that 
$\mathbb{S}-Cat^{rig}$ is not a homotopy category, it is simply a full sub-category 
of $\mathbb{S}-Cat^{\Gamma}$).
We have two functors 
$$\mathbb{S}-Cat^{rig} \longrightarrow S^{1}-SSet,$$
to the category of $S^{1}$-equivariant simplicial sets. 
The first one sends $T$ to $LT$ together with its natural action of $S^{1}$. The second
one sends $T$ to $T(1,1)$ with the trivial $S^{1}$-action. These two functors
are considered as objects in 
$Ho(Fun(\mathbb{S}-Cat^{rig},S^{1}-SSet))$, the homotopy category of
functors. Let us denote these two objects by $L : T \mapsto LT$ and 
$E : T \mapsto T(1,1)$.

\begin{conj}\label{conj2}
There exists a morphism in $Ho(Fun(\mathbb{S}-Cat^{rig},S^{1}-SSet))$
$$Tr : L \longrightarrow E,$$
in such a way that for any rigid symmetric monoidal 
$\mathbb{S}$-category $T$ the induced morphism of simplicial sets
$$Tr : LT \longrightarrow T(1,1)$$
is the trace map described above.
\end{conj}

It can be shown that conjecture \ref{conj2} implies conjecture \ref{conj1}. In fact the tensor
$2$-categories $Dg_{parf}(X)$ are the $2$-truncation of natural 
rigid symmetric monoidal $(2,\infty)$-categories, which can also be considered
as $(1,\infty)$-categories by only considering invertible higher morphisms. An application of  
the above conjecture to these rigid symmetric monoidal $(1,\infty)$-categories 
give a solution to conjecture \ref{conj1}, but this will be explained in more detailed in 
a future work. To finish this part on rigid $(1,\infty)$-categories let us mention that 
a recent work of J. Lurie and M. Hopkins on universal properties of 
$(1,\infty)$-categories of $1$-bordisms seem to solve conjecture \ref{conj2} (\cite{lu3}). We think to 
have another solution to the part of conjecture \ref{conj2} concerned with 
the rigid  symmetric monoidal $(1,\infty)$-categories of saturated dg-categories, which 
is surely enough to also imply conjecture \ref{conj1}. This again will be explained in a future work.

\begin{center} \textbf{Relations with variations of Hodge structures} \end{center}

The derived loop space $\mathrm{L}X$ and the $S^{1}$-equivariant derived category
$D^{S^{1}}_{parf}(\mathrm{L}X)$ have already been studied in \cite{bn}. In this work
the category $D^{S^{1}}_{parf}(\mathrm{L}X)$ is identified with a certain derived category 
of modules over the Rees algebra of differential operators on $X$ (when, 
say, $X$ is smooth over $k$ of characteristic zero). We do not claim to fully understand this
identification but it seems clear that objects in $D^{S^{1}}_{parf}(\mathrm{L}X)$ could be identified
with some kind of fltered complexes of $D$-modules on $X$. Using this identification
our categorical Chern character $Ch^{cat}(T)$ probably encodes the data
of the negative cyclic complex $HC^{-}(T)$ of $T$ over $X$ together with its 
\emph{Gauss-Manin connection and Hodge filtration}. In other words, 
$Ch^{cat}(T)$ seems to be nothing more than the variation of Hodge structures
induced by the family of dg-categories $T$ over $X$. As far as we know the
construction of such a structure of variations of Hodge structures on 
the family of complexes of cyclic homology associated to a family of sturated dg-categories 
is, up to conjecture \ref{conj1}, a new result (see however \cite{ge} for the construction of a Gauss-Manin connection on
cyclic homology).  We also think it is a remarkably nice fact that variations of Hodge structures appear 
naturally out of the construction of our Chern character for categorical sheaves. 

It is certainly possible to describe the cohomological Chern character of \ref{d5} using this point of view
of Hodge structure. Indeed, $HC_{0}^{-,S^{1}}(\mathrm{L}X)$
is close to be the $S^{1}$-equivariant de Rham cohomology of $\mathrm{L}X$, and using 
a localization formula it is probably possible to relate $HC_{0}^{-,S^{1}}(\mathrm{L}X)$ with 
$HP_{0}(X)[[t]][t^{-1}]$, where $HP_{0}(X)$ is periodic cyclic homology of $X$
and $t$ is a formal paramater. We expect at least a morphism
$$HC_{0}^{-,S^{1}}(\mathrm{L}X) \longrightarrow HP_{0}(X)[[t]][t^{-1}].$$
The image of $Ch^{coh}(T)$ by this map should then be closely related to 
the \emph{Hodge polynomial of $T$}, that is $\sum_{p}Ch(Gr^{p}HC^{-}(T))t^{p}$, where
$Gr^{p}HC^{-}(T)$ is the $p$-th graded piece of the Hodge filtration of Hochschild homology and
$Ch$ is the usual Chern character for sheaves on $X$. \\

\begin{center} \textbf{Back to elliptic cohomology ?} \end{center}

In the introduction of this work we mentioned that our motivation to start thinking about categorical sheaf theory 
was elliptic cohomology. However, our choice to work in the context of algebraic geometry 
drove us rather far from elliptic cohomology and it is at the moment unclear whether our work on the 
Chern character can really bring any new insight on elliptic cohomology. About this we would like to make the following 
remark. Since what we have been considering are categorified version of algebraic vector bundles, it seems
rather clear that what we have done so far might have some relations with what could be called
\emph{algebraic elliptic cohomology} (by analogy with the distinction between algebraic and topological 
$K$-theory). However, the work of M. Walker  shows that 
algebraic $K$-theory determines completely topological $K$-theory (see \cite{wa}), and that it is possible
to recover topological $K$-theory from algebraic $K$-theory by an explicit construction. Such a striking fact 
suggest the possibility that an enough well understood algebraic version of elliptic cohomology could
also provide some new insights on usual elliptic cohomology. We leave this vague idea for future works.

\end{document}